\documentclass[journal,twocolumn,10pt]{IEEEtran}

\pdfoutput=1

\usepackage{amsmath,bm}

\usepackage[]{times}
\usepackage{epsfig}
\usepackage{subfigure}
\usepackage{amsmath}
\usepackage{amssymb}
\usepackage{graphicx}
\usepackage{subfigure}
\usepackage{multirow}
\usepackage{flushend}
\usepackage{url}
\usepackage{color}
\newcommand{\bega}{\begin{eqnarray}}
\newcommand{\ega}{\end{eqnarray}}
\newcommand{\bb}{\begin{equation}}
\newcommand{\ee}{\end{equation}}

\begin{document}

\title{Reinforcing Power Grid Transmission\\
with FACTS Devices}

\author{{\bf Vladimir Frolov}$^{(1,4)}$, {\bf Scott Backhaus} $^{(2,4)}$ and {\bf Misha Chertkov} $^{(3,4,1)}$\\
$^{(1)}$ Skolkovo Institute of Science and Technology,
100 Novaya Str, Skolkovo 143025, Russia  \\
$^{(2)}$ Material Physics and Applications Division, LANL, Los Alamos, NM 87545, USA \\
$^{(3)}$ Center for Nonlinear Studies and Theoretical Division,
LANL, Los Alamos, NM 87545, USA\\
$^{(4)}$ New Mexico Consortium, Los Alamos, NM 87544, USA}

\maketitle

\begin{abstract}
We explore optimization methods for planning the placement, sizing and operations of Flexible Alternating Current Transmission System (FACTS) devices installed into the grid to relieve congestion created by load growth or fluctuations of intermittent renewable generation.  We limit our selection of FACTS devices to those that can be represented by modification of the inductance of the transmission lines. Our master optimization problem minimizes the $l_1$ norm of the FACTS-associated inductance correction subject to constraints enforcing that no line of the system exceeds its thermal limit. We develop off-line heuristics that reduce this non-convex optimization to a succession of  Linear Programs (LP) where at each step the constraints are linearized analytically around the current operating point.  The algorithm is accelerated further with a version of the cutting plane method greatly reducing the number of active constraints during the optimization, while checking feasibility of the non-active constraints post-factum. This hybrid algorithm solves a typical single-contingency problem over the MathPower Polish Grid model (3299 lines and 2746 nodes) in 40 seconds per iteration on a standard laptop---a speed up that allows the sizing and placement of a family of FACTS devices to correct a large set of anticipated contingencies.  From testing of multiple examples, we observe that our algorithm finds feasible solutions that are always sparse, i.e., FACTS devices are placed on only a few lines. The optimal FACTS are not always placed on the originally congested lines,  however typically the correction(s) is made at line(s) positioned in a relative proximity of the overload line(s).
\end{abstract}


\section{Introduction}
\label{sec:Intro}

Power grids are undergoing significant evolution, e.g. experiencing disturbances caused by fluctuating renewable generation, aging infrastructure with high costs for replacements or upgrades, and the continual growth and evolution of electrical loads.  However, the basic operational functions of the grid remains the same, e.g., to balance the generation and load in interconnected transmission grids while respecting the transmission line constraints (among several other constraints and stability limits).  These issues contribute to the growing interest \cite{97NAGA,98GH,01GCG,01OCG,03OB,03WW,05GA,06SV,06SC,07LFN,08RH,08SLSK,08WLLJ,09FMF,12SB} in utilizing a new class of smart transmission-grid devices known commonly as Flexible Alternating Current Transmission System (FACTS) devices \cite{FACTS,97FACTS-def,00HG,02MV}.  Although FACTS devices can provide fundamentally new degrees of flexibility to transmission grids, they are also expensive (especially fast solid-state devices) and must be deployed in ways that maximizes benefit to the overall system.

Congestion relief is one attractive application of FACTS devices that can provide significant system benefit for limited capital investment.  However, the non-local nature of power flows over transmission networks makes it difficult to determine where and how large a FACTS device(s) to install. To fully extract the benefit of FACTS devices, one ought to build in the control of FACTS devices into algorithms for optimal placement and sizing, i.e., formulate and solve a comprehensive FACTS planning and control problem that considers the transmission system as a whole.  This type of integration is not entirely new, see e.g. \cite{98GH,01GCG,05GA}, but a system-level approach to the problem is deemed difficult. Indeed, the system-wide optimization suffers, in general, from the ``curse of dimensionality"---the effort to exactly compute results scales exponentially with the system size and the number possible FACTS locations.  It is impractical to implement such direct approaches for large transmission systems consisting of thousands of components. However, some new approaches, brought into the power systems from statistical physics \cite{10CPS,11CSPB,12DCB} and optimization \cite{ll11opf,12BCH,low12,cvhb12approx}, suggested that even seemingly difficult optimization problems can be modified into formulations which are provably polynomial (or even linear). Such modifications are not achievable in all cases, and one often must search for computationally efficient yet empirically accurate heuristics. We take the latter path in this manuscript.

Planning/sizing and operation/control problems are normally considered separately. Planning problems seek to place FACTS devices and determine the required operational ranges. The planning problem is a strategic decision, which takes into account multiple operational conditions extending over years. On the contrary the control problem seeks the best response to current conditions or those expected within a relatively short period of time. However different the two type of problems may seem, they may be approached in a unified principled manner similar to an approached developed in in \cite{11DBC} for placement of storage. In \cite{11DBC},  storage devices where placed all nodes in the network, but with appropriate penalties or costs on utilizing the storage.  By simulating operations over an many different situations, the statistics of storage utilization were used to eliminate lightly used nodes in favor of the more frequently utilized nodes.  Repeating this process several times ``manually'' created sparse storage solutions with this external procedure.  Forcing sparsity internal to the process turns the optimization problem into a far more difficult mixed formulation with both discrete and continuous variables entering the formulation. However, recent results in compressed sensing \cite{08CWB} \footnote{See also \cite{10JC} for some related discussion in the power system context.} suggest that adjustment of the optimization cost function may allow sparsity to be enforced implicitly while keeping the optimization variables in the preferred continuous domain.

We seek to adapt some of the techniques discussed above to the problem of placing, sizing, and operating FACTS devices to relieve transmission congestion.  Here, we consider two ways of imposing congestion. The first is uniform load growth by application of the same scale factor to all loads in the transmission grid model. A simple measure of grid robustness is the amount of load scaling that can be accommodated before a transmission constraint is violated.  Alternatively, the amount load scaling achievable can be used a method to distinguish between different potential upgrades to improvement grid performance, i.e. the upgrade that enables a larger scaling is ranked higher.  There are many ways to improve the system to accommodate load growth, however traditional solutions such as building new transmission lines or generators close to load are becoming increasingly difficult---an environment that begins to favor FACTS devices with their the smaller footprint and lack of emissions.  In particular, FACTS that are able to control power flows or modify the inductance of the transmission lines are particularly attractive because a few well-placed FACTS devices can be intelligently controlled to redistribute power flows and relieve congestion on many power lines, potentially including lines remote the the FACTS devices themselves.

Throughout this manuscript, we consider uniform scaling as the measures of stress and system performance.  Other methods of imposing transmission-grid stress are fluctuations of intermittent renewable generation  and $N-1$ contingencies.  Fluctuations of intermittent generation will elicit a response from the controllable generators (via primary frequency control and subsequently via secondary control/AGC).  For large fluctuations or those with particulary unfortunate spatial configurations, the power flows driven by the response of the controllable generators may cause overloading of transmission lines.  A method for finding the set of the most probable fluctuations that cause failure have been described \cite{10CPS,11CSPB,12BH}.  Then, instead of uniform load growth, the set of instantons may be used to apply stress to the transmission grid.  $N-1$ contingencies also apply stress to the system, and congestion in potential $N-1$ networks is a typical factor leading to reduced capability of the actual system.

We model the FACTS device\cite{97FACTS-def} as a continuous modification of the inductance of a transmission line.  Such a model is directly applicable to thyristor-controlled series compensation (TCSC) devices that use a continuously controllable reactor in parallel with a bank of switchable capacitors\cite{05GA}. Phase-shifting transformers with appropriate local controls may also be considered within this model.  We seek to resolve the following questions:
\begin{itemize}
\item Can a particular infeasible configuration(s) of generation and load, i.e. a configuration that violates one or more transmission network constraints, be corrected with FACTS devices?
\item When such a correction is possible, what is the optimal (least ``expensive'') set of FACTS devices that achieves this goal?
\end{itemize}
These two questions are formulated as one non-convex optimization problem with the $l_1$ norm of the inductance deviation vector as an objective and the thermal limits for all the transmission lines as constraints.  We construct efficient heuristic algorithms that resolve these questions in a computationally-efficient manner,  and moreover, we empirically observe that the optimal solutions are always sparse, i.e. only few FACTS devices are actually needed.

Layout of material in the rest of the paper is as follows.  We pose the problem of optimal placement of the inductance correcting FACTS devices in Section \ref{sec:model}; describe solution algorithm in Section \ref{sec:algorithm}; and discuss results of our experiments in Section \ref{sec:experiments}. Section \ref{sec:conclusions} is reserved for summary of the results,  conclusions and discussing the path forward.

\section{Optimization Model}
\label{sec:model}

Before considering the optimal placement and sizing of FACTS devices, we first discuss a few preliminaries.

\subsection{General Setting: Power Flow (PF) and Optimum Power Flow (OPF)}
\label{subsec:setting}

The transmission grid is represented as a graph ${\cal G}=({\cal V},{\cal E})$ consisting of a set of vertices ${\cal V}$ and a set of undirected edges ${\cal E}$.  The graph is assumed fixed. The vector of transmission line susceptances is  $\beta\doteq(\beta_{ab}|\{a,b\}\in{\cal E})$, and it is this $\beta$ that we assume is modified by FACTS devices and is the subject of optimization.  The transmission lines have power flow limits represented by the vector $f\doteq(f_{ab}=f_{ba}>0|\{a,b\}\in{\cal E})$.  The generators have over and under generation limits $\underline{g}\doteq(\underline{g}_a|a\in{\cal V})$, $\overline{g}\doteq(\overline{g}_a|a\in{\cal V})$. We split the set of nodes ${\cal V}={\cal V}_c \cup {\cal V}_u$ into controllable nodes (traditional generators participating in control or perhaps modern controllable loads) and  uncontrollable nodes (traditional loads and renewable generators). In our notation, the number of vertexes  and edges is  $N=|{\cal V}|$ and $M=|{\cal E}|$, respectively.

The transmission grid power flows are computed using the DC approximation, i.e., (a) the resistance of transmission lines can be ignored in comparison with the line inductance, (b) the normalized voltage magnitude is unity at all nodes, and (c) phase difference between neighboring nodes is small, although we will use an iterative approximation to relax this approximation.  The relationship between phases $\theta_a$ and powers $p_a$ consumed/injected at any node $a$ becomes linear:
\begin{eqnarray}
&& \forall a\in{\cal V}:\quad p_a =\sum_{b\sim a}\beta_{ab}(\theta_a-\theta_b)=(L\theta)_a,
\label{DC}\\ &&
L_{ab}\doteq\left\{
\begin{array}{cc}
\sum_{c\sim a} \beta_{ac},& a=b;\\
-\beta_{ab},& a\neq b,\ \{a,b\}\in{\cal E};\\
0, & \mbox{otherwise}.
\end{array}
\right.
\label{L}
\end{eqnarray}
where, $L\doteq(L_{ab}=L_{ba}|a,b\in{\cal E})$ is the $N\times N$ weighted graph-Laplacian matrix.  The total balance of the real power is $\sum_{a\in{\cal V}} p_a=0$, and the phase vector $\theta=(\theta_a|a\in{\cal V})$ is defined up to a global shift.

In the rest of this manuscript, we find it convenient to express $L$  as $L(\beta)=\nabla^+*\beta*\nabla$. Here, $\nabla$ is the $(2M)*N$ rectangular matrix which generates a $2M$ vector of phase differences across the undirected edges $(a,b)$, i.e., $\sum_{c=1}^N \nabla_{(a,b);c}\theta_c=\theta_a-\theta_b$.   Also, $\beta$ is a $(2M)\times (2M)$ diagonal matrix with elements $\beta_{ab}$ for both directed edges $(a,b)$ and $(b,a)$.  Finally, $\nabla^+$ is the $N\times (2M)$ transpose of $\nabla$.

The phase vector is then expressed as $\theta=L^{-1}p$ where $L^{-1}$ is the pseudo-inverse of $L$.  The pseudo-inverse is well-defined because, by construction, $N-1$ eigenvalues of $L$ are positive with the final zero eigenvalue corresponding a uniform shift of $\theta$.  This ambiguity in the phase vector $\theta$ is customarily resolved by fixing the phase at the slack node (often the largest generator in the system).

Our objective in this manuscript is to optimally use FACTS devices to relieve transmission congestion caused by additional stress applied to a base grid configuration either via load growth or fluctuating generation or load. To ensure that our results are not biased by a poor choice of the base configuration, we select this configuration by performing an optimal power flow (OPF) calculation.  Starting from a  known or forecasted set of uncontrolled power injections $p_u=(p_a|a\in{\cal V}_u)$, we solve the following DC-OPF problem:
\begin{eqnarray}
&\min\limits_{p_c}& \sum_{a\in{\cal V}_c} c_a(p_a);\label{DC-OPF}\\
&\mbox{s.t.}& \forall \{a,b\}\in{\cal E}:\  \left|((\beta*\nabla)*L^{-1}*p)_{(a,b)}\right|\leq f_{ab},\nonumber\\
&& \forall a\in{\cal V}_c:\  \underline{g}_a\leq p_a\leq\overline{g}_a,\nonumber
\end{eqnarray}
where the cost of generation $c_a(p_a)$ is a non-decreasing function of $p_a$.  In the expression inside the absolute value on the left and side of the first set of constraints, $p$ is understood to be the $N\times 1$ column vector of the $p_a$, and this expression reduces to $\beta_{ab}(\theta_a-\theta_b)$, i.e. the DC approximation of the power flow on transmission line $(a,b)$.  The slightly more complicated form in Eq.~\ref{DC-OPF} will prove useful in the discussion of optimization over the FACTS devices.
When considering stress applied by uniform load growth, a DC-OPF could be run to re-optimize the controllable generation at each level of load growth.  Instead, we uniformly grow all generation by the same scale factor that we use to grow loads. This way we overload some lines and then we try to find the solution (fix the lines) by searching for optimal susceptances with our algorithm.

\subsection{Stressing of the Optimal Solution}
\label{subsec:stress}

Assuming a feasible solution  $p^*$ to the DC-OPF exists, one way to measure the robustness of the solution/grid is to scale all of the generation and loads  by a constant factor $\alpha>1$ and test to see if the solution still respects the constraints in Eq.~\ref{DC-OPF}, i.e., if
\begin{eqnarray}
\forall \{a,b\}\in{\cal E}:&\quad \alpha \left|((\beta*\nabla)*L^{-1}*p^*)_{\{a,b\}}\right|\leq f_{ab}. \label{alpha-edge}
\end{eqnarray}
The critical $\alpha$, i.e $\alpha_c$, is the smallest value of $\alpha$ where at least one of the transmission line constraints in Eqs.~(\ref{alpha-edge}) is violated.  A super-critical, infeasible grid configuration is then given by $p^{(in)}=\alpha p^*$ with $\alpha>\alpha_c$.  Note that we assume the system is only constrained by transmission congestion and that generation constraints do not play a significant role.

A second way to stress the solution $p^*$ to achieve an infeasible $p^{(in)}$ is to choose certain configurations
\cite{10CPS,11CSPB} of the uncontrolled fluctuating generation and loads $p_u$ that, when combined with the response of the controlled generators $p_c$ participating in secondary frequency control, results in violations of transmission line constraints.  Much like the uniform scaling discussed above, a one-dimensional family of $p^{(in)}$ can be generated by scaling the original $p_u$ fluctuation.

\subsection{Problem Formulation}
\label{subsec:formulation}

We seek to improve the grid robustness to stress (e.g., increasing the value of $\alpha_c$) by optimally placing FACTS devices. Stress can be applied to the base grid configuration $p^*$ in many ways to create an infeasible $p^{(in)}$, including the two discussed above.  However, the problem formulation that follows is independent how the stress applied.

In this initial work, we only consider FACTS devices whose effect can be approximated as modifying the susceptance (inductance) of a transmission line.  Such FACTS devices are expensive, and we incorporate this via the cost function $H(\beta;\beta^{(in)})$,  where $\beta^{(in)}$ and $\beta$ are the vectors of susceptance for the infeasible solution and for the solution after FACTS modifications. We model $H(\beta;\beta^{(in)})$ as the $l_1$-norm of the FACTS-induced change in $\beta$ and pose the following optimization problem:
\begin{eqnarray}
&&\hspace{-1cm}\min\limits_{\beta}\sum_{\{a,b\}\in{\cal E}}|\beta_{ab}-\beta_{ab}^{(0)}|
\label{placement}\\
&&\hspace{-1cm}\mbox{s.t. } \forall \{a,b\}\in{\cal E}:\ \left|((\beta*\nabla)*(L(\beta))^{-1}*p^{(in)})_{ab}\right|\leq f_{ab}.
\nonumber
\end{eqnarray}
In spite of some useful features of the dependence of the graph-Laplacian on $\beta$, see e.g. \cite{08GBS}, the nonlinear conditions in Eq.~(\ref{placement}) are generally non-convex in $\beta$. (See Appendix \ref{subsec:three} where the non-convexity is illustrated on a three node example.) We resolve this complication  in Section \ref{sec:algorithm} by designing a greedy but efficient heuristic algorithm that enables solution of (\ref{placement}) over large practical instances.

Eq.~(\ref{placement}) is not guaranteed to have a feasible solution, i.e. FACTS corrections of line inductance may not be sufficient to correct all constraint violations when the system is severely stresses.  It is interesting to discover the amount of stress to reach this second critical boundary.  For case of uniform scaling in Section \ref{subsec:stress}, we seek to find the $\alpha_c^\prime>\alpha_c$ such that $p^{(in)}= \alpha_c^\prime p^*$ results in infeasibility of Eq~\ref{placement}.  This analysis will be discussed on examples in Section \ref{sec:experiments}.

Finally, we note that is does not make sense to apply (\ref{placement}) to radial distribution grids because, when ${\cal G}$ is a tree, fixing $p$ leads to an unambiguous, $\beta$-independent set of power flows. However,  even in the simplest case of a single loop the optimization problem (\ref{placement}) becomes interesting and nontrivial (see discussion of Appendix \ref{subsec:three}).

\subsection{Critical Discussion}
\label{subsec:critical}

Before discussing implementation details in Section \ref{sec:algorithm}), we first comment on the assumptions and limitations, but also advantages, natural generalizations and uses of the FACTS placement problem stated in Eq.~(\ref{placement}).  The following basic assumptions/caveats were used in the formulation of Eq.~(\ref{placement})
\begin{itemize}
\item[(1)] The power flow equations are linearized, i.e. we use the DC approximation.
\item[(2)] The $\ell_1-$norm cost function in Eq.~(\ref{placement}) is over-simplified as it does not include fixed costs of installing FACTS that do not scale with $\beta-\beta^{(in)}$.
\item[(3)] Sparsity of FACTS placement is a desirable property of any solution of Eq.~(\ref{placement}), but sparsity is not explicitly encouraged by the formulation.
\item[(4)] Generation constraints $\underline{g}_a$ and $\overline{g}_a$ are ignored.
\item[(5)] Only a single configuration $p^{(in)}$ is considered as a guidance for FACTS installation.
\item[(6)] The FACTS-corrected network is not guaranteed to be $N-1$ reliable (or reliable with respect to a more general list of contingencies).
\end{itemize}

We comment of these items in order.

(1) Linearizing the power flow equations is advantageous as it allows us to express the transmission-line thermal constraints with an analytic formula in terms of $\beta$. Generalization to the exact nonlinear power flows is straightforward,  however in this case getting an explicit formula for the line flows in terms of $\beta$ may be a challenge.  The lack of an explicit formula complicates the efficient computation of the $\beta$-derivative and iterative linearization over $\beta$ used to derive an efficient algorithm in Section \ref{sec:algorithm}.  A viable alternative for the exact nonlinear power flow equations is iterative linearization over {\it both} susceptance and phase angles around a base-solution.

(2) The choice of the $l_1$ norm is illustrative.  Many other cost functions may be incorporated without creating significant complications, including those with graph-inhomogeneity that encourage building new FACTS at specific locations.

(3) Although we not to explicitly encourage sparsity, our experiments with large networks (see Section \ref{sec:experiments}) naturally result in sparse solutions.  We conjecture that the $l_1$ norm may promote sparsity in ways similar those in compressive sensing \cite{08CWB}.

(4) Generalization of Eq.~(\ref{placement}) accounting for generation limits is straightforward.

(5) We have developed efficient algorithms for solving Eq.~\ref{placement} so that we may simultaneously consider many different stressed configurations by simply replicating Eq.~\ref{placement}.  Specifically, we create $n$ replicas of the constraints in Eq.~(\ref{placement}) for $n$ different stressed configurations, $p^{(1)},\cdots p^{(n)}$, replacing $M$ constraints in Eq.~(\ref{placement}) by $n*M$ conditions and minimizing the cost function in Eq.~(\ref{placement}) over the extended set of constraints.

(6) $N-1$ reliability can be incorporated using methods similar to (5) above.  Specifically, for a ``virtual grid'' with a one line disconnected (but the other $M-1$ lines still functioning), we derive the set of $M-1$ constraints that are explicitly dependent on all but one component of $\beta$.  We combine all the constraints into a set of $(M-1)*M$ constraints which, additional to $M$ base constraints, must be satisfied by a solution of Eq.~(\ref{placement}).

The optimization problem (\ref{placement}), and its generalizations briefly discussed above, can be used both for planning of new FACTS installations into existing grid and for controlling existing FACTS devices. In the case of planning, the list of contingencies $p^{(1)},\cdots p^{(n)}$ could include a number of different configurations of a future stressed grid (including $N-1$ contingencies) for both the system peak and minimum.  The solution of (\ref{placement}), if it exists, would then provide the lowest cost placement and sizing of FACTS that would resolve all of the network constraint violations for the stressed configurations considered. In the case of operations and control, contingencies may be considered one by one, e.g. a single $N-1$ contingency or one of the worst-case instantons discussed in \cite{10CPS,11CSPB}.  Equation~(\ref{placement}) then enables a reactive control that quickly computes the FACTS set points to maintain a feasible solution.

\section{Optimization Algorithm}
\label{sec:algorithm}

Eq.~(\ref{placement}) presents two challenges: the inequality constraints are nonlinear and (b) there are too many constraints -- thousands for a large transmission grid like that discussed in Section \ref{sec:experiments}. In Sections \ref{subsec:linear} and \ref{subsec:cutting}, we will describe ways overcome these challenges independently, and in Section \ref{subsec:synthesis} present a synthesis strategy combines the two approaches into one iterative procedure.

\subsection{Linearization of Constraints}
\label{subsec:linear}

We take advantage of the explicit expression of the conditions in Eq.~(\ref{placement}) and linearize these around the current value of $\beta^*$,
\begin{eqnarray}
&& \forall \{a,b\}\in{\cal E}:\ ((\beta*\nabla)*(L(\beta))^{-1}*p^{(1)})_{ab}\nonumber\\ &&
\approx  A_{ab}(\beta^*)+\sum_{\{c,d\}\in {\cal E}}(\beta_{cd}-\beta_{cd}^*)B_{ab;cd}(\beta^*).
\label{linear}
\end{eqnarray}
The M-dimensional vector $A$ and $(M\times M)$ matrix $B$ can be pre-computed explicitly (or tabulated in more general case) for a range of values of $\beta^*$. The resulting Sequential Linear Programming (SLP) \cite{SLP} approximate iterative algorithm then becomes\\

\underline{\bf Direct Algorithm}
\begin{itemize}
\item {\bf Start:} $k=0$. Initiate with $\beta^{(0)}=\beta^{(in)}$.
\item {\bf Step 1:} Linearize conditions in Eq.~(\ref{placement}) about $\beta^{(k)}$ according to Eq.~(\ref{linear}) and solve the following LP, $\beta^{(k+1)}=$
\begin{eqnarray}
&&\hspace{-1cm}\mbox{argmin}_\beta\sum_{\{a,b\}\in{\cal E}}|\beta_{ab}-\beta_{ab}^{(0)}|
\label{LP-linear}\\
&&\mbox{s.t. } \forall \{a,b\}\in{\cal E}:\nonumber\\
&&\Biggl|A_{ab}(\beta^{(k)})\!+\!\sum_{\{c,d\}\in {\cal E}}\!(\beta_{cd}\!-\!\beta^{(k)}_{cd}) B_{ab;cd}(\beta^{(k)})
\Biggr|\leq f_{ab},
\nonumber
\end{eqnarray}
\item {\bf Step 2:} If $|\beta^{(k+1)}-\beta^{(k)}|$ is larger than a convergence tolerance, then set $k=k+1$ and go to Step 1.
\item {\bf End:} Output $\beta^{(k+1)}$ as the solution.
\end{itemize}
Note that in actual experiments discussed in Section \ref{sec:experiments}, the algorithm converges after a small number of iterations with $\beta^{(k+1)}$ exactly equal to $\beta^{(k)}$.

\subsection{Cutting Plane}
\label{subsec:cutting}

One problem with the optimization in (\ref{placement}) and the LP reformulation in Eq.~(\ref{LP-linear}) is the complexity due to the large number of constraints. However, in practical cases, very few of these constraints are actually violated, suggesting that the complexity of the brute force implementation can be drastically reduced with standard cutting-plane algorithms, see e.g. \cite{99Ber,03Avd}.  A modification of Eq.~(\ref{LP-linear}) consists of cutting the constraints in two groups, ``included" and ``excluded", and updating the two groups till convergence is reached. Initially, the included group consists of the constraints which are violated for $\beta=\beta^{(k)}$, while all other constraints are placed in the excluded group. At every step of the inner iteration (with respect to the algorithm explained in Section \ref{subsec:linear}), we solve Eq.~(\ref{LP-linear}) using only the included constraints, and then check if any of the excluded constraints are violated in the resulting solution.  If no excluded constraints are newly violated, we conclude that an optimal solution of the full problem is found.  Otherwise, we pick the most violated constraint and move it from the excluded group to the included group and repeat. This simple and straightforward procedure works very effectively in all the practical cases we tested, normally stopping in less than a handful of steps.

\subsection{Synthesis}
\label{subsec:synthesis}

Our numerical experiments suggest that the algorithms performance can be boosted significantly by alternating the outer-loop LP solution with the inner-loop cutting plane, instead of waiting for convergence of the inner-loop cutting plane step:\\
\underline{\bf Improved Algorithm}
\begin{itemize}
\item {\bf Start:} $k=0$. Initiate with $\beta^{(0)}=\beta^{(in)}$. Split the list of $2M$ inequality constraints in Eq.~(\ref{placement})in to the included (${\cal E}^{(in)}$) and excluded (${\cal E}^{(out)}={\cal E}\setminus {\cal E}^{(in)}$) groups.
    \item  {\bf Step 1a:} Linearize the ${\cal E}^{(in)}$ Eq.~(\ref{placement}) about $\beta^{(k)}$ and solve the following LP, $\beta^{(k+1)}=$
\begin{eqnarray}
&&\hspace{-1cm}\mbox{argmin}_\beta\sum_{\{a,b\}\in{\cal E}}|\beta_{ab}-\beta_{ab}^{(0)}|
\label{LP-synt}\\
&&\mbox{s.t. } \forall (a,b)\in{\cal E}^{(in)}:\quad f_{ab}\geq\nonumber\\
&& \sigma_{ab}\Biggl(A_{ab}(\beta^{(k)})\!+\!\sum_{\{c,d\}\in {\cal E}}\!(\beta_{cd}\!-\!\beta^{(k)}_{cd}) B_{ab;cd}(\beta^{(k)})
\Biggr),
\nonumber
\end{eqnarray}
where $\sigma_{ab}$ is $+1$ or $-1$ depending on the signature of the respective directed (one-sided) inequality.
\item  {\bf Step 1b:} Update ${\cal E}^{(in)}$ moving all the currently violated constraints that are in ${\cal E}^{(out)}$  for $\beta^{(k)}$ to ${\cal E}^{(in)}$.
\item {\bf Step 2:} If $|\beta^{(k)}-\beta^{(k)}|$ is larger than tolerance
or if the update set on the previous step was not empty, then set $k=k+1$ and
go to Step 1a.
\item {\bf End:} Output $\beta^{(k+1)}$ as the solution.
\end{itemize}

\begin{figure}[h!]
\centering
\includegraphics[width=0.5\textwidth]{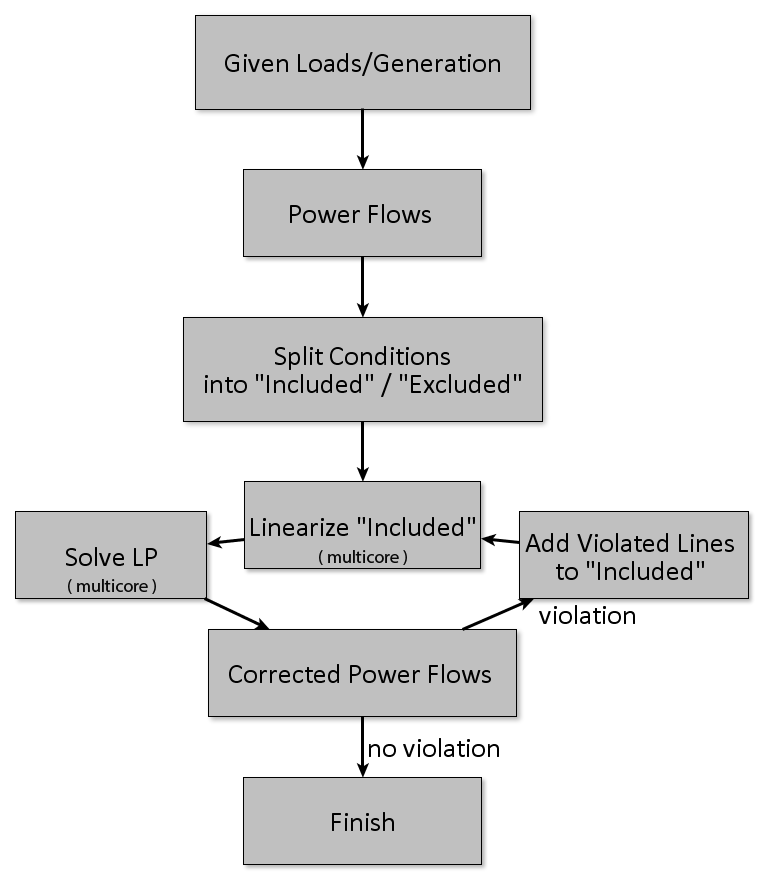}
\caption{Flowchart of the improved algorithm. Multicore procedures are marked. The most time consuming procedure is "Solve LP".
\label{fig:flow_chart}}
\end{figure}

It is straightforward to verify that the fixed point of this improved algorithm procedure will also be a fixed point of the direct algorithm procedure (wait till convergence of the inner loop before making the next outer loop step). Given the global non-convexity of Eq.~(\ref{placement}) potential landscape, one obviously cannot guarantee that the fixed point of the improved procedure will always coincide with a fixed point of the direct procedure. However, the improved procedure is as good as the direct one for finding a local minimum --- which is exactly the problem we are aiming to solve.

\section{Experiments}
\label{sec:experiments}

Next, we illustrate the algorithms of Section \ref{sec:algorithm} on two transmission grids.  The first is a 30-bus model from the software Matpower \cite{matpower} \footnote{To convert MatPower cases into the standard format (matrix L and vector p), transformers and phase shifters are turned off, double lines are combined to form one line with value of throughput and inductance calculated from two lines. Double generators also were combined to form one generator.} (see Fig.~\ref{fig:30_overload1}), a small enough grid that we can develop some intuition about how the FACTS are being utilized.  The second is the Polish grid where we consider two versions -- a 2737-bus summer grid and a 2746-bus winter model (also available in Matpower, see Fig.~\ref{fig:Polish1}).  Tests of the improved algorithm of Section \ref{subsec:synthesis} revealed that the number of iterations required for convergence is unexpectedly small -- less than a dozen for all the cases we experimented with (analysis of convergence is shown in the Fig.2). In the case of the Polish grid, each iteration of our  algorithm took ~30 seconds on a standard quad-core processor.

\begin{figure}[h!]
\centering
\includegraphics[width=0.5\textwidth]{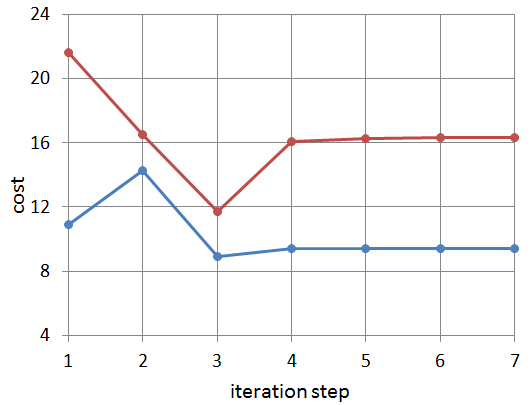}
\caption{Convergence of the cost with iteration step for 30-bus model (blue graph) and for the Polish model (red graph).
\label{fig:30_overload1}}
\end{figure}

The FACTS placements selected by the solutions of the numerical experiments also revealed a number of interesting features of the solutions that we discuss below.

\subsection{Non-locality}
\label{subsec:nonlocal}

The nonlocal behavior is illustrated on two examples of the 30-node model with $\alpha/\alpha_c=1.4$ ($40\%$ overload) in Fig.~\ref{fig:30_overload1} and $\alpha/\alpha_c=1.9$ ($90\%$ of overload) in Fig.~\ref{fig:30_overload2}.  Initially overloaded lines are marked in red and modified lines are marked in green. The numbers near the modified lines are percentage of the susceptance change.  Very often, although not always, our algorithm chose NOT to decrease the susceptance of the overloaded lines to restrict the power flow on them, but instead it modifies the susceptance of nearby lines to reroute power flows around the congested transmission lines.  This rather frequent nonlocal behavior of the optimal solutions suggests that the optimal placement of FACTS devices is nontrivial problem and calls for a computationally efficient approach, such as the one we consider here, so that the method can be  applied to much larger grids as well.

Arrowheads on the lines in Fig.~\ref{fig:30_overload2} indicate the direction the original $\alpha/\alpha_c=1.9$ power flows and the smaller arrows (green/down or red/up) indicate whether the original power flows decreased or increased after the FACTS were placed.  The original power flows out of generator G1 overloaded line L1.  Our algorithm drives $\beta\rightarrow 0$ on line L4 and increases $\beta$ on line L2 effectively rerouting the power from G1 to the lower right of the network.  The increase in power flow on L2 also relieves the overload on L3 demonstrating the inherently non-local effects of FACTS.  This redistribution of power flows has even longer range effects.  The major reduction of flow on line L4 has two other beneficial impacts.  First, it draws more power from G2 (in spite of the decrease in $\beta$ on L5) relieving the overload on L6.  In addition, it forces a reversal of the power flow on line L7 relieving the overload on L8.

The small test grid in Fig.~\ref{fig:30_overload2} enables us to build up some intuition about the non-local effects that affect the optimal placement and sizing of FACTS.  We next consider the same uniform load scaling for the much larger Polish grid.  Figure~\ref{fig:Polish1} highlights the small region of the Polish grid where all of the overloads and susceptance modifications occur for $\alpha/\alpha_c = 1.38$. The results are presented in the tabular form for different values of $\alpha/\alpha_c$ in Fig.~\ref{table:Polish_overload} and demonstrate behavior similar to the much smaller 30-bus network. For small $\alpha/\alpha_c$ up to at least 1.04, the optimal solution is local, i.e. our algorithm selects to simply reduce the susceptance of line 375, which reduces the flow on this overloaded line.  As $\alpha/\alpha_c$ grows, non-local behavior becomes apparent.  For $\alpha/\alpha_c \ge 1.1$, additional lines become overloaded, however, none of these lines are selected for susceptance modification.  Interestingly, line 375, which was the initial overloaded line, is selected for susceptance modification in all of the solutions.

The details of the highlighted region of Fig.~\ref{fig:Polish1} are shown in Fig.~\ref{fig:Polish1a} for $\alpha/\alpha_c=1.38$.
The structure of the solution is similar to the 30-bus configuration one. The overloads are resolved by changing the susceptances of the near lines to reorganize power flows structure in the system. It can be seen that the corrections are not always non-local. Here the 375th line was overloaded initially and also corrected in the solution.

We note that in some of the cases, our algorithm selected susceptance corrections that set a line's total susceptance to zero, i.e. corresponding to the removal of the line from the network. This curious fact was also verified directly by manually removing the line in question and rerunning our algorithm.  The resulting solution was indeed a valid solution (i.e. no lines overloaded). The structure of such solution is the same. If we do not consider line 33 which was removed (shown in the Fig.4), the same lines were changed and the corrections are approximately the same.

\begin{figure}[h!]
\centering
\includegraphics[width=0.5\textwidth]{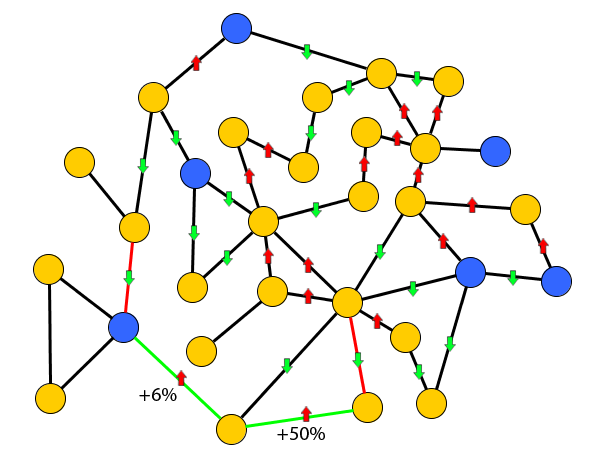}
\caption{Visualization of the 30-node model illustrating the non-locality of the optimal solution with $\alpha/\alpha_c=1.40$ yellow buses - loads, blue buses - generators. Lines marked in red were overloaded in the base case (two lines), while those marked in green where selected by our algorithm for susceptance modification by FACTS devices (two lines). The percentage of the correction to the susceptance is shown next to the corrected lines. The difference in power flows through the lines after the corrections are shown with short green arrows (decrease of the flow) and short red arrows (increase of the flow). If there is no arrow then power flow is the same.
\label{fig:30_overload1}}
\end{figure}

\begin{figure}[h!]
\centering
\includegraphics[width=0.5\textwidth]{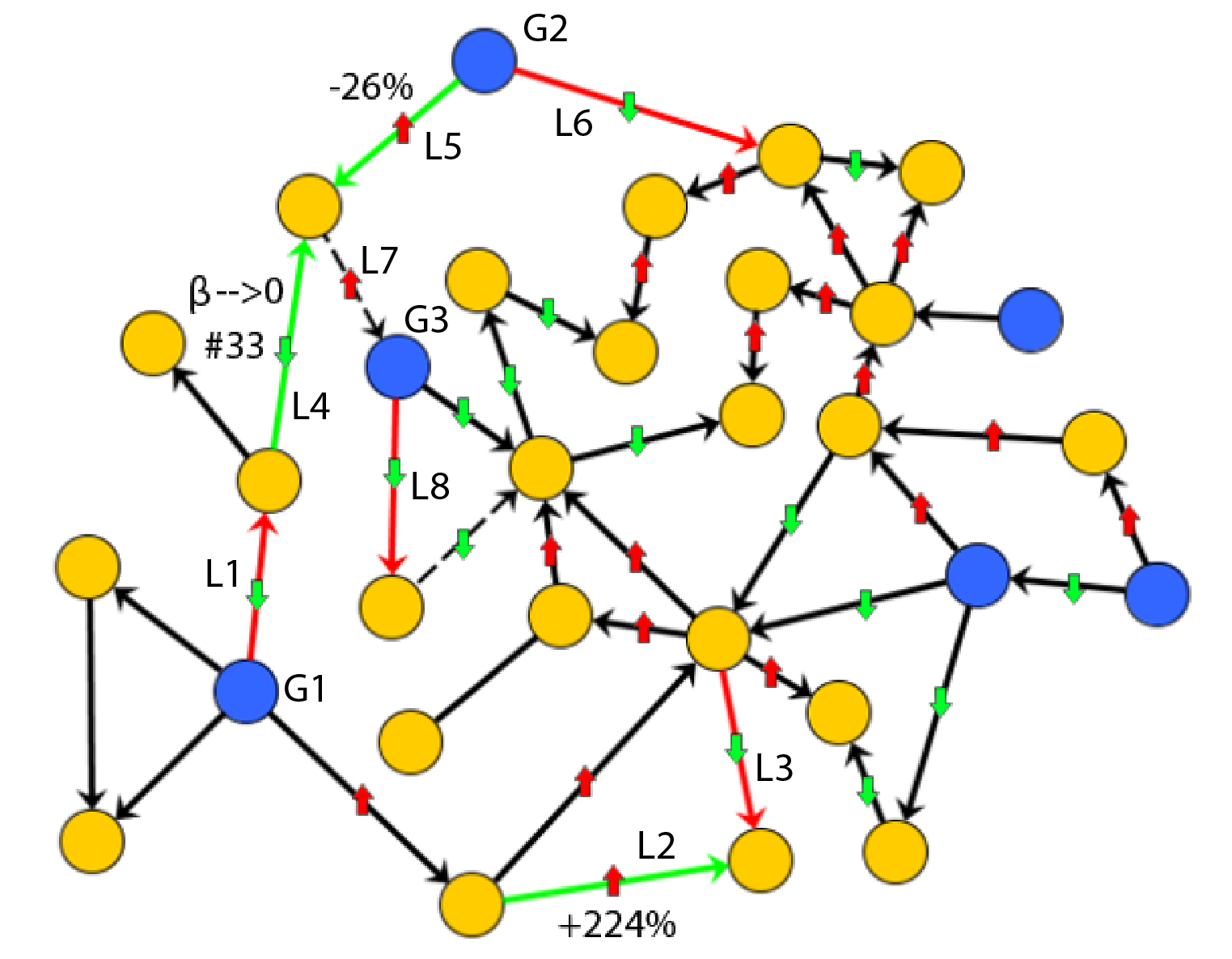}
\caption{Visualization of the 30-node model illustrating the non-locality of the optimal solution with $\alpha/\alpha_c=1.90$ yellow buses - loads, blue buses - generators. Lines marked in red were overloaded in the base case (four lines), while those marked in green where selected by our algorithm for susceptance modification by FACTS devices (three lines). Directions of the in-line arrows indicate directions of the flow prior to the FACTS corrections. The percentage of the correction to the susceptance is shown next to the corrected lines (susceptance of the 33th line is going to 0). The difference in power flows through the lines after the corrections are shown with green short arrows (decrease of the absolute value of the flow) and red short arrows (increase of the absolute value of the flow). If there is no arrow then power flow is the same. If power flow through the line changes the direction after installation of FACTS (corrections made) it is illustrated by dotted line (two lines).
\label{fig:30_overload2}}
\end{figure}

\begin{figure}[h!]
\centering
\includegraphics[width=0.5\textwidth]{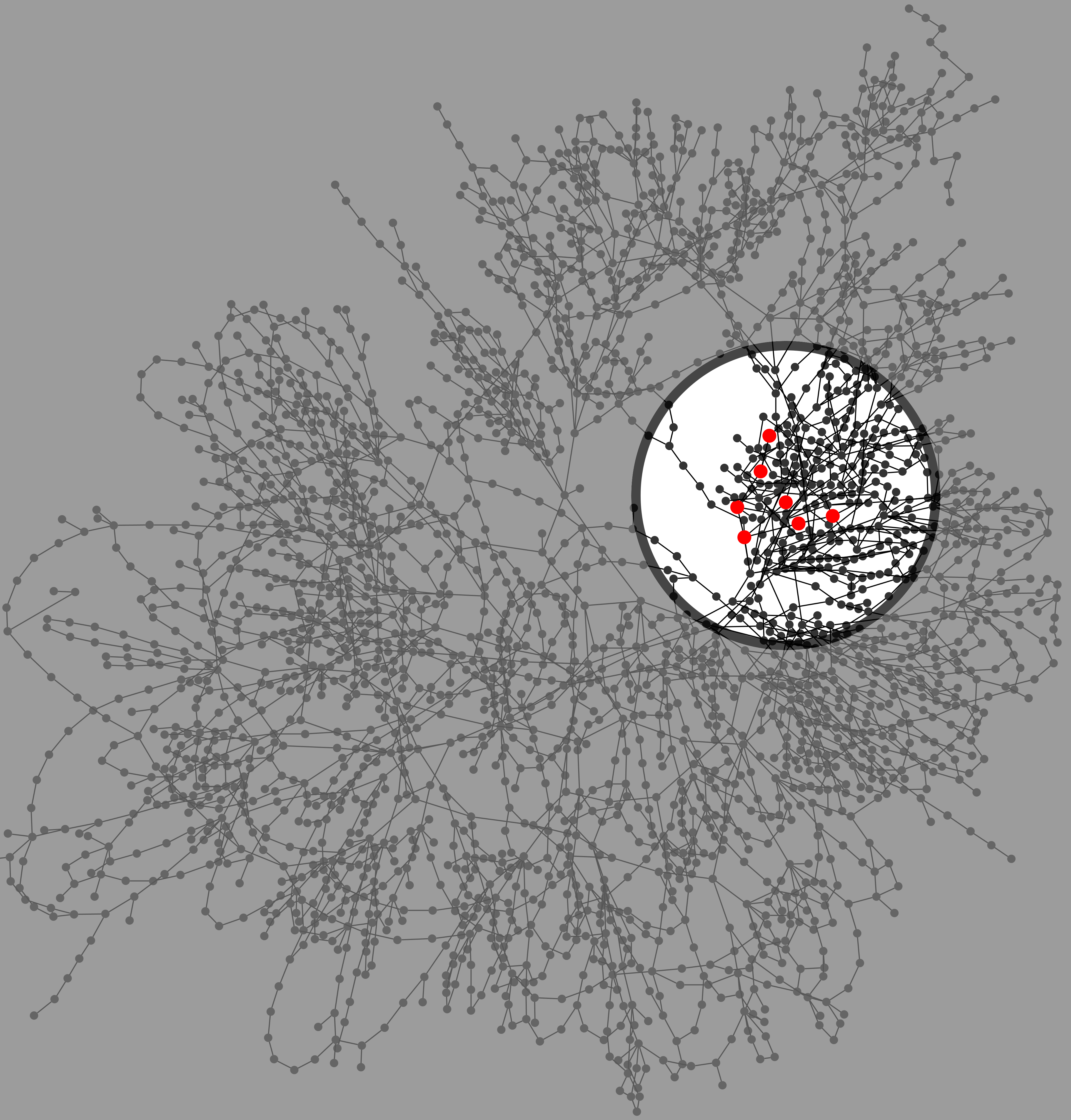}
\caption{Non-geographical visualization of the entire Polish grid.  The highlighted region contains all of the overloaded and FACTS-modified lines for the cases shown in the table of Fig.~\ref{table:Polish_overload}.  Part of the highlighted region is shown in more detail in Fig. \ref{fig:Polish1a} for $\alpha/\alpha_c=1.38$ with the red nodes here corresponding to the big nodes of Fig.~\ref{fig:Polish1a}.
\label{fig:Polish1}}
\end{figure}

\begin{figure}[h!]
\centering
\includegraphics[width=0.5\textwidth]{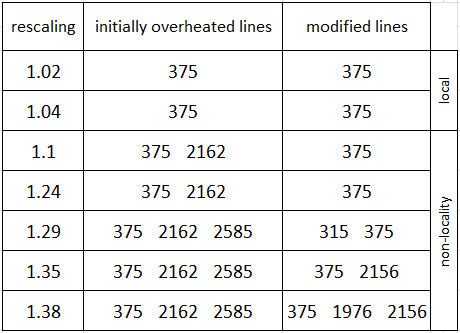}
\caption{Table showing which lines are overloaded in the base solution of the Polish grid in Fig.~\ref{fig:Polish1}  for selected values of $\alpha/\alpha_c$ and lines selected by our algorithm for susceptance modification by FACTS.  For $\alpha/\alpha_c\le 1.04$, the overload is relieved by simply lowering the susceptance on the overloaded line to limit the flow on that line.  For larger $\alpha/\alpha_c$, the solution become nonlocal with the additional susceptance modifications placed on lines that are not overloaded in the base solution.
\label{table:Polish_overload}}
\end{figure}

\begin{figure}[h!]
\centering
\includegraphics[width=0.5\textwidth]{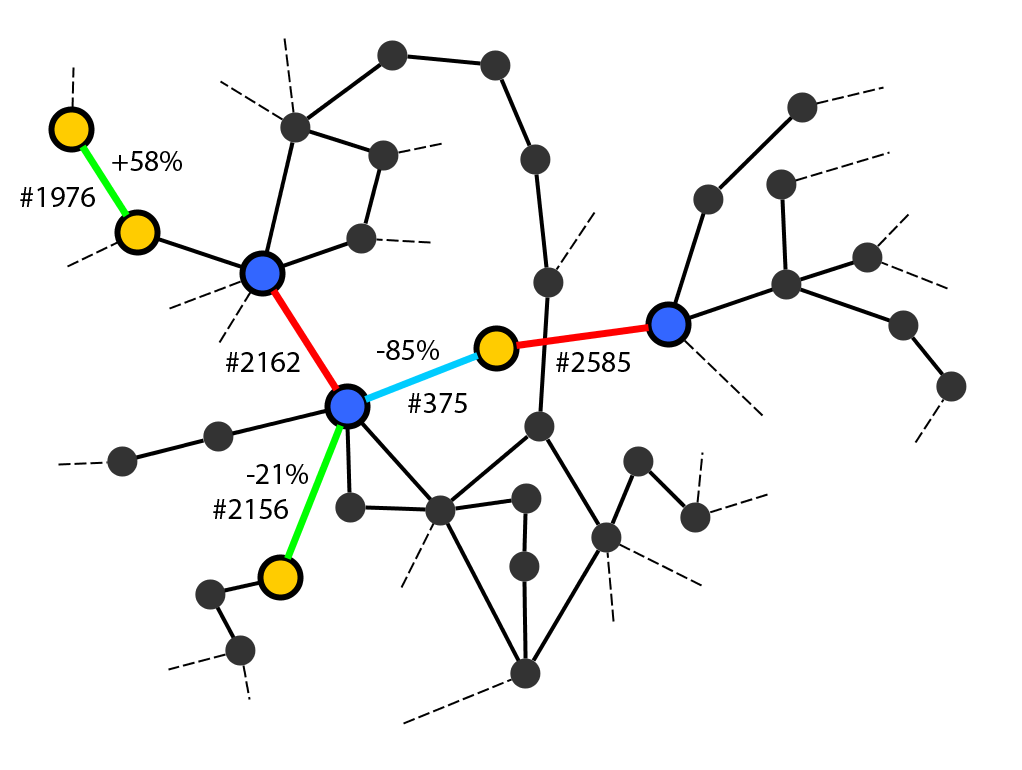}
\caption{Detail of the highlighted region of the Polish grid from Fig.~\ref{fig:Polish1} showing the lines overloaded in the base case ($\alpha=$1.38) in red and those selected by our algorithm for susceptance modification in green. Generators are marked as blue and loads are as yellow nodes. The one line shown in blue was both overloaded and modified.  The numbers labelling these lines correspond to the line numbers in the table of Fig.~\ref{table:Polish_overload}. Although the Polish grid is much larger, the behavior is very similar to that of the 30-node grid in Figs.~\ref{fig:30_overload1} and \ref{fig:30_overload2}, i.e. susceptances are primarily modified on lines that are near to the originally overload lines to encourage the power flow to be rerouted around the overloads to relieve congestion.  \label{fig:Polish1a}}
\end{figure}

\subsection{Sparsity}
\label{subsec:sparsity}

A second interesting observation concerns sparsity of the optimal solution.  Instead of requiring small modifications of many lines throughout the network, our algorithm makes significant susceptance modifications to a only few lines with the number of modified lines typically the same or slightly smaller than the number of overloads in the base case. This behavior is qualitatively similar in  examples of the 30-node network and of the Polish networks and can be seen in Figs.~\ref{fig:30_overload1}, \ref{fig:30_overload2}, and \ref{fig:Polish1a}.  We note that our $\ell_1-$norm cost function does not explicitly promote sparsity, i.e. for a given budget of susceptance modification, it does not cost more to spread it out over the entire network rather than concentrate it on a few lines.  However, the sparsity of the solution emerges naturally.  We conjecture that this sparsity is a natural property of the ``N-1 redundancy'' engineered into electrical networks, i.e. N-1 redundancy generally requires that there be at least two paths to deliver power to loads, and if one of paths becomes overloaded, an increase in susceptance of an alternate path will deliver more power thus relieving the congestion on the first path \footnote{Even thought the sparsity was not enforced directly, the use of the $l_1$ norm in the cost function of Eq.~(\ref{placement}) was most probably helping implicitly, similarly to how the sparsity arises in compressed sensing, see e.g. \cite{}.}.  These observations suggest an additional cutting plane-like heuristic that could speed up our algorithm even further: instead of optimizing over the susceptances of all of the lines, one could restrict our attention to the set of lines that are near to the overloaded lines.

\subsection{Uniform Load Scaling and Emergence of Local Optima}
\label{subsec:rescaling}

By efficiently solving the optimization problem in Eq.~\ref{placement}, our algorithm allows a more thorough exploration of the solution space.  In particular, we investigate the emergence of multiple minima of the cost function in Eq.~\ref{placement} as we uniformly re-scale the loads by $\alpha$ in both the 30-node and Polish networks.

\begin{figure}[h!]
\centering
\includegraphics[width=0.5\textwidth]{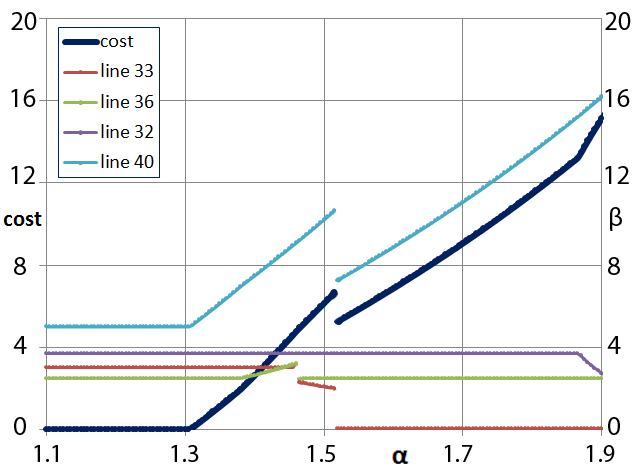}
\caption{ Dependence of the optimal cost on the scaling factor, $\alpha$, shown for the working case of the 30 node model. Lines which were corrected and their final susceptances are illustrated for lines $\#33$, 36, 32 and 40.
\label{fig:cost_30}}
\end{figure}

\begin{figure}[h!]
\centering
\includegraphics[width=0.5\textwidth]{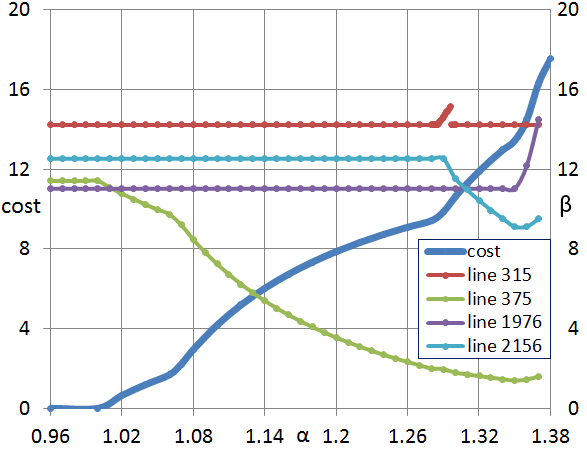}
\caption{ Dependence of the optimal cost on the scaling factor, $\alpha$, shown for one of the summer cases of the Polish model. Lines which were corrected and their final susceptances are illustrated. Bends of the cost graph correspond to the corrections of additional lines.
\label{fig:cost_Polish}}
\end{figure}

For the 30-node network, Fig.~\ref{fig:cost_30} clearly shows a jump in the cost at $\alpha = 1.51$ signaling emergence of multiple (at least two, but possibly more) local minima.  Investigation of the details of the solutions above and below the jump shows that the corrections to susceptances of lines $\#33$ and $\#40$  are significantly different.  Above the jump, the net susceptance of line $\#33$ becomes approximately zero indicating that this line has effectively been removed from the network. The results in Fig.~\ref{fig:cost_30} leads us to suspect that, by always initiating with the base solution, our algorithm has become trapped in a local minima for $\alpha<1.51$.  To verify this suspicion,  we initiated the algorithm with a configuration equivalent to the bare configuration for all lines but with line $\#33$ removed.  Dependence of optimal cost of $\alpha$ for this case is illustrated in the Fig.10 with final susceptances of the lines shown. Assuming that turning the line off costs zero we significantly decreased the cost of corrections. But the jump of the cost still exists which means that we still have another minima of the cost function.

To prove the fact that there are more than one minima (more than one final solutions for one initial configuration of the susceptances) let us compare the initial solution (Fig.8) and the solution with line $\#33$ turned off (now we will calculate the cost of such operation in terms of the cost function defined). Costs of two different solutions for the same initial grid (30-bus grid) are illustrated in the Fig.11.

Three regions can be seen in the graph. In the 3rd region, the solutions are the same. In the 1st region, there is extra cost associated with the forced reduction of the susceptance of line $\#33$.  This effective removal of line $\#33$ does not result in lower cost.   But in the 2nd region it can be seen that this removal does decrease the final cost compared to a solution where line $\#33$ was not removed. There are two different solutions with different costs implying that at leats two minima of cost function exist.

\begin{figure}[h!]
\centering
\includegraphics[width=0.5\textwidth]{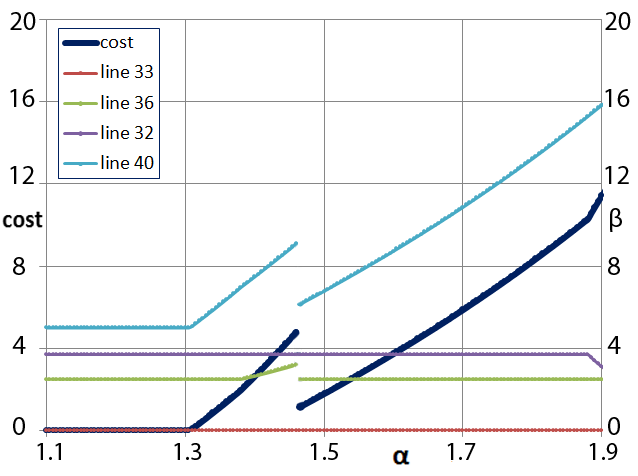}
\caption{ Dependence of the optimal cost on the scaling factor, $\alpha$, shown for the working case of the 30 node model. Lines which were corrected and their final susceptances are illustrated for lines $\#33$, 36, 32 and 40. Line $\#33$ was turned off with zero cost.
\label{fig:lines_Polish}}
\end{figure}

\begin{figure}[h!]
\centering
\includegraphics[width=0.5\textwidth]{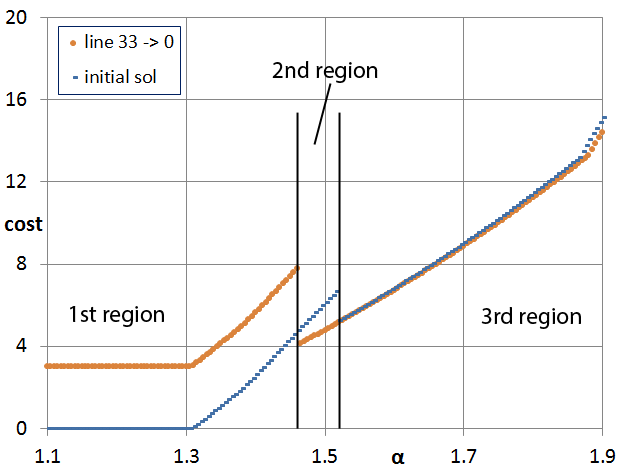}
\caption{ Comparison of the solutions given in the Fig. 8 and Fig. 10. Solutions are evaluated in terms of the cost function defined above. Turning line $\#33$ off is not free now.
\label{fig:cost compare}}
\end{figure}

The same analysis (but no jumps in this case) for Polish case is shown in the Fig.9. Different bends of the cost vs scaling factor dependence seen in Fig.~\ref{fig:cost_Polish} also indicates competition between different optimal structures.  The structures are different in the number of overloaded lines, the number FACTS-corrected lines (Fig.6), and the magnitudes of the corrections. Indeed, at small $\alpha$ the optimal configuration contains one overloaded and one corrected line.  At $\alpha\approx 1.06$, another line becomes overloaded, however the configuration is still correctable with a single FACTS device positioned at the same line as before. At $\alpha\approx 1.26$, another line becomes overheated, but it does not require an additional FACTS device until $\alpha\approx 1.29$, etc.

\subsection{Robust Optimal Placement: Correcting Multiple Configurations}
\label{subsec:multiple}

Up to this point, we have only considered using FACTS to modify line susceptances to correct a single configuration that causes network overloads.  However, as discussed in Section~\ref{subsec:critical}, we can easily robustly optimize the placement and sizing of FACTS to correct  $n$ configurations that lead to overloads.  We modify the LP portion of the improved algorithm of Section \ref{subsec:synthesis}  by extending the list of constraints. For each of the $n$ network configurations, $p^{(1)},\cdots p^{(n)}$, on every iteration step (Fig.1) we create a list of $2m^{(1)},\cdots 2m^{(n)}$ inequalities for the $m^{(i)}$ transmission lines from ``included" and combine all these inequalities in one extended list. We then replace the list of the directed edge labels ${\cal E}$ in Eq.~\ref{placement} with the new composite list and iterate the improved algorithm as described in Section~\ref{sec:algorithm}.

To demonstrate the effectiveness of this method, we consider a simple situation where we simultaneously optimize FACTS placement and sizing over two supercritical cases.  The first case describes uniformly scaled Polish grid's winter off-peak configuration from MatPower. The second case is generated from the preceding configuration by adding a fraction $X$ to the previous loads where $X$ is distributed uniformly from $-0.7$ to $0.7$.  Following the load changes, the second configuration includes a generation adjustment via an optimal power flow, and then the loads were  uniformly scaled.

\begin{figure}[h!]
\centering
\includegraphics[width=0.5\textwidth]{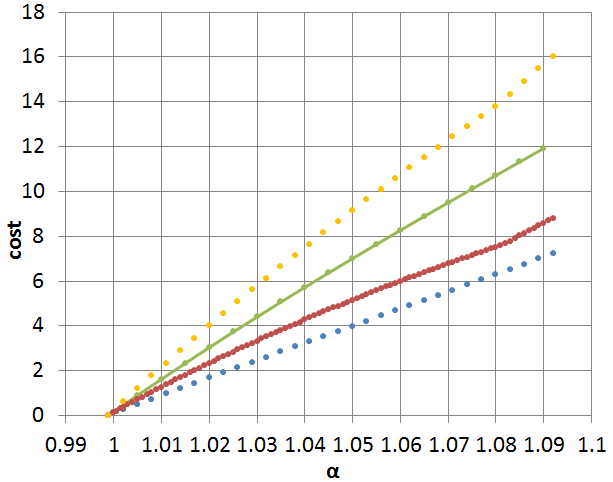}
\caption{Optimal cost vs scaling parameter shown for FACTS placements over three different settings of Polish grid are marked in red, blue and green colors. The result of optimization over a base, off-peak winter configuration is shown in red. The blue curve shows the optimization result for another configuration with the same line topology as in the base case,  however with randomly generated corrections to loads and respective OPF adjustment of the generation. The green curve shows robust-optimal FACTS placement which optimal of all placements which remain feasible for both,  base and corrected, configurations.  This green curve should be compared with naive combination of the two single-configuration optimization shown in yellow.
\label{fig:joint_scaling1}}
\end{figure}

Fig.~\ref{fig:joint_scaling1} shows the resulting cost of the robust optimization (green curve)  compared to the optimal FACTS placement and sizing that only considers one or the other scenario (red and blue curves).  Since the robust optimization is correcting both overloads at the same time, the total cost is higher than either the red or blue curves.  However, it is lower than the yellow curve which is the sum of the maximum cost per line over the two scenarios, i.e. the cost one would naively compute by optimizing FACTS placement and sizing for the two scenarios independently.

The advantage of the robust optimization is not observed in all cases.  For example, if in the two scenarios the overloaded lines were spatially well separated, the lack of interaction between the FACTS devices and transmission line constraints at these locations would effectively split the robust optimization back into two separate, uncoupled optimizations.  A second case where the robust optimization does not improve the results is when the scenarios  have exactly the same set of line overloads.  However, if the lines that are overloaded in the different scenarios are in reasonable proximity, the FACTS devices that correct one scenario may be leveraged to correct another scenario at lower total cost.  For the results displayed in Fig.~\ref{fig:joint_scaling1}, we have purposely picked a randomly-generated second scenario that displays some overlap with the original, uniformly-scaled ($\alpha=1.08$) Polish winter base case.

\section{Conclusions and Future Works}
\label{sec:conclusions}

In this manuscript we studied how adjusting power line inductances with FACTS devices can help to reduce congestion in the power systems.  Our main finding is the suggestion of efficient heuristics which can be used for both operational (short term) and planning (long term) purposes. The heuristics are used in an algorithm that minimizes the deployment of FACTS using an $l_1$-norm as a measure FACTS cost.  The minimization is done under the condition that thermal overloads are relieved for a specially selected bad configuration(s) of load and generation. The algorithm performance is illustrated on moderate and large scale models. We find that the resulting optimal FACTS deployments are sparse, i.e. necessary at only a handful of lines.  However,  the lines to be corrected are not necessarily lines with the most sever overload (without the correction applied).

There are number of directions for future work in this area including:
\begin{itemize}
\item Incorporating in the operational scheme multiple dangerous configurations and possibly including an outer step to search for these configurations.
\item Considering FACTS control in combination with other devices and controls, e.g. energy storage and generation re-dispatch.
\item Extension of the planning paradigm to stochastic setting that takes into account the statistics of multiple stressed configurations and considers the investment problem of placing and sizing FACTS devices to minimize risks.
\item Generalize the model presented here to  account for more general (and accurate) AC modeling of power flows,  thus modeling risks of loss of synchrony and voltage collapse (in addition to the risk of the thermal overload discussed in this manuscript).
\end{itemize}

\section*{Appendix}

\subsection{Motivating Three Node Example}
\label{subsec:three}

\begin{figure}[h!]
\centering
\includegraphics[width=0.3\textwidth]{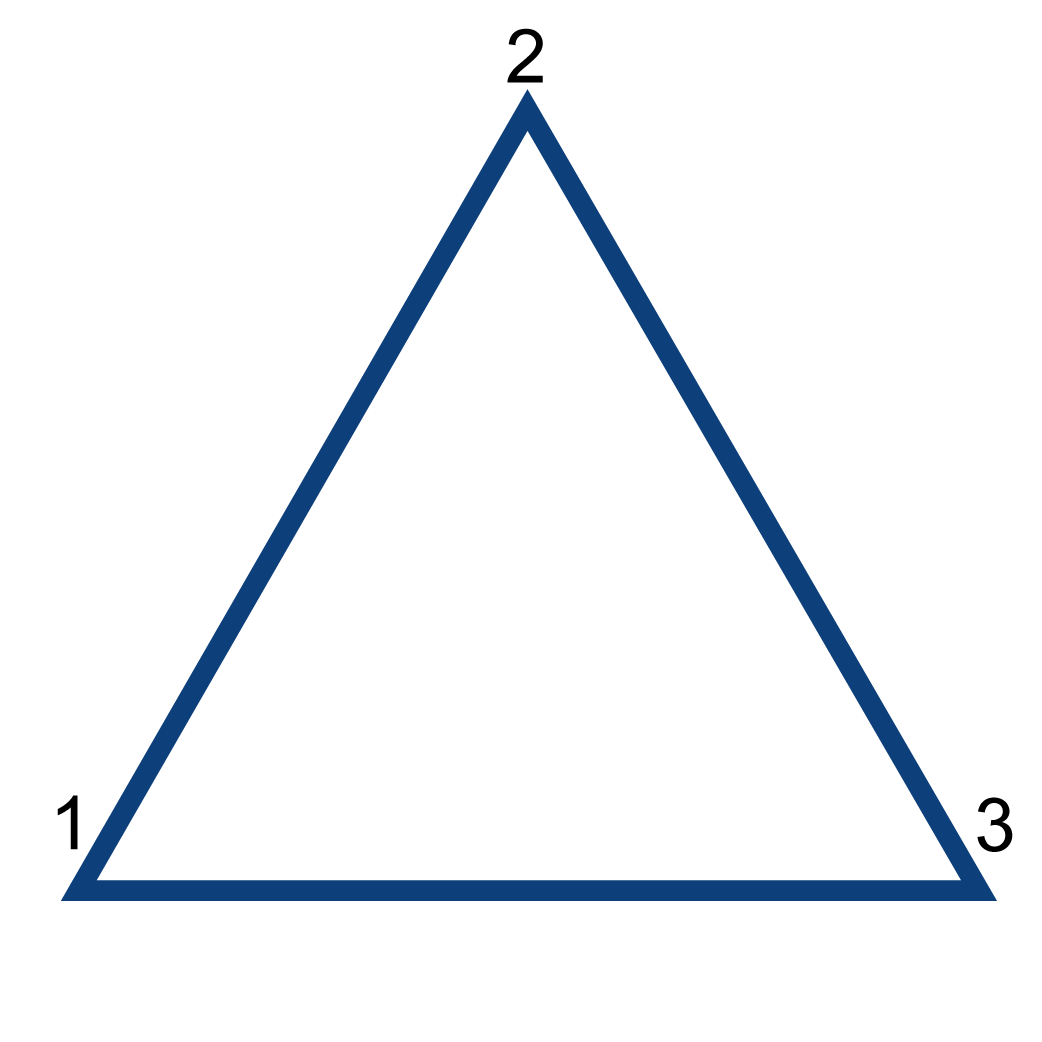}
\centering
\caption{The simple graph of the three node example.}
\label{fig:3}
\end{figure}

\begin{figure}[h!]
\centering
\includegraphics[width=0.45\textwidth]{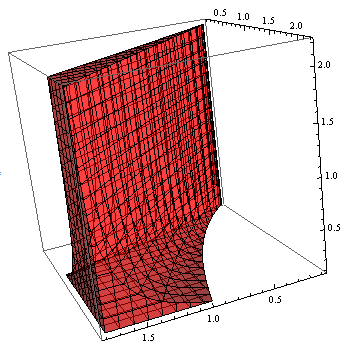}
\centering
\caption{Visualization of the three dimensional $\beta$ space for the three node example of Fig.~\ref{fig:3}.  Non-convexity of the domain is clearly seen.}
\label{fig:3_1}
\end{figure}

\begin{figure}[h!]
\centering
\includegraphics[width=0.45\textwidth]{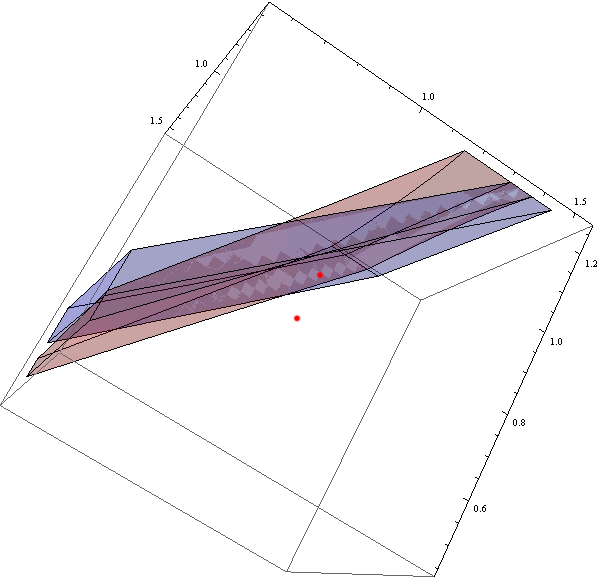}
\centering
\caption{The three dimensional $\beta$ space as in Fig.~\ref{fig:3_1} but rotated. Red area marks feasible domain of the original nonlinear conditions in (\ref{placement}) for the three node example,  and blue area corresponds to its linearized version as in Eq.~(\ref{LP-linear}). Red dots stand for the initial, $\beta_0$, configuration (seen outside of the domains),  and the final solution (inside of the domains) achieved in one iteration of the direct algorithm of Section \ref{subsec:linear}.}
\label{fig:3_2}
\end{figure}

Here we illustrate non-convexity of our general formulation (\ref{placement}) on the simple three node example of Fig.~\ref{fig:3}, where $a=1,\cdots,3$.
Power flows are solved explicitly and then the domain limited by six conditions in Eq.~(\ref{placement}) is shown, as the function of three susceptances, $\beta_{12},\beta_{13},\beta_{23}$, in Fig.~\ref{fig:3_1}, where we see clearly nonlinearity and non-convexity of the domain in $\beta$. Fig.~\ref{fig:3_2} illustrates our linearized algorithm (no need to add cutting plane trick here).  Here, the exact optimization domain (equivalent to the one shown in Fig.~\ref{fig:3_1} but rotated for better view) is shown in red, and domain linearized around the basic $\beta_0$ case is shown in blue. Two red points mark initial state (outside of the domains) and final optimal state (inside of the domains) found in one iteration. 

\subsection{Sensitivity of Line Flow to Local Change of Susceptance}

It is instructive to analyze sensitivities of the line flows to changes of the line suspectances on  this simple three node example. Power flow over a line of the triangle, say line $(1,2)$, is
\begin{equation}
p_{12}=\frac{\beta_{12}(p_1\beta_{23}-p_2\beta_{13})}{\beta_{12}\beta_{13}+\beta_{12}\beta_{23}+\beta_{13}\beta_{23}},
\label{p12}
\end{equation}
where $p_i$ is the production/consumption at the node $i$ (which is positive/negative). Then the sensitivity of the flow to the change of susceptance at the same or neighboring lines (under fixed nodal productions/consumptions) becomes
\begin{eqnarray}
&& \frac{\partial p_{12}}{\partial \beta_{12}}=
\frac{\beta_{13}\beta_{23}p_{12}}{\beta_{12}\left(\beta_{12}\beta_{13}+\beta_{12}\beta_{23}+\beta_{13}\beta_{23}\right)}
\label{sens_12_12}\\
&& \frac{\partial p_{12}}{\partial \beta_{23}}=
\frac{\beta_{12}\beta_{13}p_{23}}{\beta_{23}\left(\beta_{12}\beta_{13}+\beta_{12}\beta_{23}+\beta_{13}\beta_{23}\right)}.
\label{sens_12_23}
\end{eqnarray}
From these formulas (and also taking into account that all the suscpetances are positive) one concludes that in order to reduce the value of the flow over a line, in the case when correction of only a single line susceptance is allowed, one either (a) decreases susceptance of the same line which is overloaded, or (b) increases/decreases susceptance of the neighboring line depending on if the directions of the original and neighboring line flows towards their common point are the same/different. (The ``same" and ``different" are interpreted as $\rightarrow\cdot\leftarrow$ and $\rightarrow\cdot\rightarrow$, respectively.)

\section*{Acknowledgment}

The work at LANL was carried out under the auspices of the National Nuclear Security Administration of the U.S. Department of Energy at Los Alamos National Laboratory under Contract No. DE-AC52-06NA25396.
This material is based upon work partially supported by the National Science Foundation award \# 1128501, EECS Collaborative Research ``Power Grid Spectroscopy" under NMC.

\bibliographystyle{IEEEtran}
\bibliography{Bib/convex,Bib/FACTS,Bib/robust,Bib/new}

\end{document}